\documentclass[11pt,a4paper,twoside]{article}
\usepackage[utf8]{inputenc}
\usepackage[T1,T2A]{fontenc}
\usepackage[russian,english]{babel}
\usepackage{graphicx}
\usepackage{array}
\usepackage[doi]{samgtu-bib}
\usepackage{mathtools}
\mathtoolsset{showonlyrefs}
\usepackage{desclist}
\usepackage{euscript}
\usepackage{mathrsfs}
\usepackage{multicol}
\usepackage{mathrsfs}
\usepackage{cite}
\usepackage{qrcode}
\allowdisplaybreaks[4]
\usepackage[nodayofweek]{datetime}

\graphicspath{{./00PIC/}}
\usepackage{wrapfig}

\begin{document}





{\bf
Delta-problems for the gene\-ra\-lized  Euler--Darboux equation\\
}


{\bf
I.~N.~Rodionova, V.~M.~Dolgopolov, M.~V.~Dolgopolov\\
}
keywords:
generalized Euler–Darboux equation, boundary value problem\\


Abstract. Degenerate hyperbolic equations are dealing with many important issues for applied nature. While a variety of degenerate equations and boundary condi\-ti\-ons, successfully matched to these differential equation, most in the characteristic coordinates reduced to Euler--Darboux one.
Some boundary value problems, in particular Cauchy problem, for the specified equation demanded the~introduction of special classes in which formulae are simple and can be used to meet the new challenges, including Delta-problems in squares that contain singularity line for equation coefficients with data on adjacent or parallel sides of the square.
In this short communication the generalized Euler--Darboux equation with negative parameters in the rectangular region is considered.

\newdate{fa}{14}{07}{2017}
\date{\displaydate{fa}}
\newdate{fb}{08}{09}{2017}
\newdate{fc}{18}{09}{2017}

\newdate{fe}{22}{09}{2017}


\section
{Introduction.\\ Statement of the problem and boundary conditions}
Degenerate hyperbolic equations occur in many important problems of dynamical systems and in questions of applied nature: the theory of infinitesimal bending of surfaces of revolution, the membrane theory of shells, in plasma magnetohydrodynamics, gas dynamics.
With all the variety of degenerated equations and boundary condi\-ti\-ons it is successfully matched to the given differential equation, the latter equation in the characteristic coordinates reduces to the Euler--Darboux equations.

Some boundary value problems (Cauchy problem, in particular) for the specifi\-ed equations  require the
introduction of special classes in which the formula for the solution becomes more simple in form and can be used to solve new tasks, including Delta($\Delta$)-problems in the squares containing the singularity line of the equation coefficients with the data on adjacent or parallel sides of the square (directed by A.~M.~Nakhushev).
%
%

The first works on Delta($\Delta$)-problems on sets representing the union of two characteristic triangles of hyperbolic equations were works of T.~S.~Kal\-me\-nov~\cite{dolg:1}, V.~F.~Volkodavov, A.~A.~Andreev~\cite{dolg:2},
A.~M.~Nakhushev~\cite{dolg:3}.
The development of Delta-problems was done in the number of works of other authors, from which we should mention~\cite{dolg:5, dolg:6,  dolg:7, dolg:8}.

Unlike the previous in the present work the formulation of problems $\Delta_2$ are on the set, which includes four of the characteristic triangles.
The~generalized Euler--Darboux equation with negative parameters is considered:
\begin{equation}
U_{\xi\eta} - \frac{p}{\eta- (\mathop{\rm sgn} \eta) \cdot \xi}  U_\xi +
\frac{p}{(\mathop{\rm sgn} \eta) \cdot \eta - \xi} U_\eta - (\mathop{\rm sgn} \eta) \cdot \lambda U =0,
\label{dolg:eq1}
\end{equation}

\begin{wrapfigure}{O}{0.45\textwidth}
\centering
\includegraphics[scale=.5]{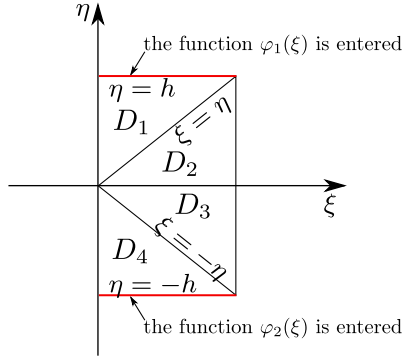}
\caption*{Region $D$}
\end{wrapfigure}

\noindent
$0<p<{1}/{2}$, $|\lambda|<\infty$ in the rectangular region $D$ bounded by characteristics of equa\-tion~\eqref{dolg:eq1} $\xi=0$, $\xi=h$, $\eta=h$, $\eta=-h$ ($h>0$), containing within itself two lines of singularity of coefficients $\eta=\xi$ and $\eta=-\xi$.
For equation \eqref{dolg:eq1}  in the region $D$ the formulation of boundary value problems $\Delta_2$ with the given values of the sought solution on the parallel sides of the rectangle, with the conditions of conjugation with respect to the solution and its normal derivative as lines of singularity of the coefficients, and the internal characteristic of the line are studied. The unique solvability of problems is proved by integral equations method. Problems are solved in the special class of functions introduced by authors in~\cite{dolg:9}.

\smallskip

\section{The solution of the problem}
Considered boundary conditions for equation (1) in the rectangular region $D$ at the parallel sides of the rectangle are of the following form:
$$ U(\xi, h) = \varphi_1(\xi), \quad U(\xi, - h) = \varphi_2(\xi), \qquad 0\leq \xi \leq h.$$
On lines of the singularity of $\eta=\xi$, $\eta=-\xi$ and on characteristic $\eta=0$ the continuity of the sought solution is gluing on.

Relatively to normal derivatives two cases are considered of the pairing on the lines $\xi=\pm\eta$. In the first case Frankl's condition of occlusion  is imposed (Problem~$\Delta_2^*$):
$$ \nu_1(\xi) = \lim_{\eta \to \xi+0} (\eta-\xi)^{-2p} (U_\xi - U_\eta) =  - \lim_{\eta \to \xi-0} (\xi-\eta)^{-2p} (U_\xi - U_\eta) = - \nu_2(\xi), $$
$$ \nu_3(\xi) = \lim_{ - \eta \to \xi-0} (\eta+\xi)^{-2p} (U_\xi + U_\eta) =  - \lim_{ - \eta \to \xi+0} (-\xi-\eta)^{-2p} (U_\xi + U_\eta) = - \nu_4(\xi), $$
$$ ~~ 0<\xi< h. $$

In the second case the gluing is carried out on the continuity of normal derivatives $\nu_1=\nu_2$; $\nu_3=\nu_4$ (Problem $\Delta_2$).
In both tasks, $\Delta_2$ and $\Delta_2^*$, in characteristic $\eta=0$ pairing sets:
$$
\lim_{\eta \to 0 + 0} \Bigl(\frac{\partial U}{\partial \eta} -  \frac{\partial U}{\partial \xi}\Bigr) =
\lim_{\eta \to 0 - 0} \Bigl(\frac{\partial U}{\partial \eta} +  \frac{\partial U}{\partial \xi}\Bigr).
$$

The basis for the decision of tasks in view is taken, obtained by the authors~\cite{dolg:9, dolg:10, dolg:11}, the solution of Cauchy problem of the special class $R_h$, which is in one of four characteristic triangles that make up the region $D$, has the form ($0<\xi<\eta<h$):
\begin{multline}
U(\xi, \eta) = \int _\eta^h T_1(s) (s-\xi)^p (s-\eta)^p   \,\,
F_{\hspace{-4mm}0 \hspace{3mm}1} \bigl(1+p; \lambda (s-\xi) (s-\eta)\bigr) ds +
\\
 + \int^\eta_\xi N_1(s) (\eta-s)^p (s-\xi)^p   \,\,
F_{\hspace{-4mm}0 \hspace{3mm}1} \bigl(1+p; - \lambda (\eta-s) (s-\xi)\bigr) ds,
\end{multline}
where $N_1= k_1T_1 - k_2 \nu_1$; $k_1$, $k_2 = \rm const$; $$\,\,
F_{\hspace{-4mm}0 \hspace{3mm}1} \bigl(\alpha; z\bigr)= \sum_{n=0}^\infty \frac{z^n}{(\alpha)_n n!}.$$

Formulas for Cauchy problem solutions in three other characteristic triangles are not given. They also contain an unknown functions $T_k$, $N_k$, $k= 2, 3, 4$, which are searched in the class of continuous in the interval $(0, h)$ and absolutely integrable functions on $[0, h]$.

The solution of the Problem $\Delta_2^*$ is reduced to the set of integral equations of the form
$$ \int _\xi^h T(s) (\xi-s)^p  \,\,
F_{\hspace{-4mm}0 \hspace{3mm}1} \bigl(1+p; \lambda s (s-\xi)\bigr) ds = \Phi(\xi, \lambda), $$
the unique solvability of which takes place when the following conditions are imposed on the given functions:
$$ \varphi_i(\xi) \in C^{(2)}
 [0, h], \qquad \varphi_i(0) = \varphi_i^\prime(0) = 0, \qquad i = 1, \, 2, $$
$$ \varphi_i(\xi) = (h-\xi)^{1+p+\varepsilon} \varphi_i^\ast(\xi), \qquad \varepsilon>0, \qquad \int_0^h \varphi_i(s) (h-s)^{-p-2} ds=0. $$

When you run
these conditions the only solution of the problem $\Delta_2^*$ is given explicitly.

The complete study of the problem of $\Delta_2$ has managed to get only if $\lambda=0$. In this case, its solution is reduced to a set of integral equations of the first kind with Cauchy kernel:
$$ \int _0^h \frac{\mu_i(s) ds}{s-\xi} = \Phi_i^\ast(\xi), \quad i= 1, \, 2;  \qquad \mu_1 = (h-s)^p s^p [\nu_1- \nu_3],$$
\begin{equation}
\label{dolg:eq2}
\mu_2 = (h-s)^p s^{p-1} [\nu_1 + \nu_3].
\end{equation}
${\Phi_i}^\ast(\xi)$ depends only on the given functions $\varphi_i$.

\smallskip

\section{Discussions of solutions}
Following the theory of singular integral equa\-tions~\cite{dolg:12}, conditions imposed on the given functions $\varphi$ under which there is the solution of equations \eqref{dolg:eq2} (not unique), and solvability conditions of equations \eqref{dolg:eq2}, and, consequently, of the Problem~$\Delta_2$.

\begin{itemize}
\item[1.] If $\varphi_i(\xi)\in C[0, 1]$, $\varphi_i^{\prime \prime} \in C(0, 1)$, $\varphi_i^{\prime \prime}$ are absolutely integrable on $[0, h]$, ${i=1, 2}$, then solutions of equations \eqref{dolg:eq2}  have the form~\cite{dolg:12}:
$$ \mu_1 = - \frac{1}{\pi} \sqrt{\frac{\xi}{h-\xi}}
\int _0^h  \sqrt{\frac{h-y}{y}} \cdot \frac{\Phi_1 ^\ast(y) }{y-\xi} dy,
$$
$$ \mu_2 = - \frac{1}{\sqrt{\xi (h-\xi)}} \biggl(\frac{1}{\pi}  \int  _0^h  \sqrt{(h-y) y} \cdot \frac{\Phi_2^\ast(y) }{y-\xi} dy + A_0\biggr), \quad A_0 = \rm const. $$

\item[2.] If the existing conditions to add
$$ \varphi_1^\prime(h) + \varphi_2^\prime(h)  = 0, \quad \mbox{but} ~ ~ \varphi_i^\prime(h)  \not= 0, \quad i = 1, \, 2,$$
then the second of equations \eqref{dolg:eq2} is uniquely solvable:
$$ \mu_2 = -  \frac{1}{\pi} \sqrt{\frac{h-\xi}{\xi}} \int  _0^h  \sqrt{\frac{y}{h-y}} \cdot \frac{\Phi_2^\ast(y) }{y-\xi} dy. $$

\end{itemize}

In conclusion, we note that for the problem definition authors were inspired by results, published in the work of I.~V.~Volovich, O.~V.~Groshev, N.~A.~Gusev, E.~A.~Kuryanovich~\cite{dolg:13}.


\EngRef



$\;$

\vspace{-10mm}


\vspace{-7mm}


\renewcommand{\baselinestretch}{0.9}

\end{document}